\documentclass[10pt, a4paper]{article}

\usepackage{amsfonts, latexsym, amsmath, amssymb, fancyhdr, amsthm, graphicx, subfigure}
\usepackage[english]{babel}

\newtheorem{theorem}{Theorem}[section]
\newtheorem{lemma}{Lemma}[section]
\newtheorem{proposition}{Proposition}[section]

\title{\textbf{A Mathematical Study of the Hematopoiesis Process with
Applications to Chronic Myelogenous Leukemia}}
\author{Mostafa Adimy$^{\dagger}$, \quad Fabien Crauste$^{\dagger}$ \quad and \quad Shigui Ruan\footnote{\small Research was partially supported by the NSERC of Canada
and the College of Arts and Sciences at the University of Miami.}
}
\date{Year 2004}

\begin{document}

\maketitle

\begin{center}
{\large $^{\dagger}$}\emph{Laboratoire de Math\'ematiques Appliqu\'ees, FRE 2570,} \\
\emph{Universit\'e de Pau et des Pays de l'Adour,}\\
\emph{Avenue de l'universit\'e, 64000 Pau, France.}\\
\emph{E-mail: mostafa.adimy@univ-pau.fr, fabien.crauste@univ-pau.fr}\\
\quad\\
{\large $^*$}\emph{Department of Mathematics, University of
Miami,}\\\emph{P. O. Box 249085, Coral Gables, FL 33124-4250,
USA.}\\\emph{ E-mail: ruan@math.miami.edu}
\end{center}

\quad

\begin{abstract}
This paper is devoted to the analysis of a mathematical model of
blood cells production in the bone marrow (hematopoiesis). The
model is a system of two age-structured partial differential
equations. Integrating these equations over the age, we obtain a
system of two nonlinear differential equations with distributed
time delay corresponding to the cell cycle duration. This system
describes the evolution of the total cell populations. By
constructing a Lyapunov functional, it is shown that the trivial
equilibrium is globally asymptotically stable if it is the only
equilibrium. It is also shown that the nontrivial equilibrium, the
most biologically meaningful one, can become unstable via a Hopf
bifurcation. Numerical simulations are carried out to illustrate
the analytical results. The study maybe helpful in understanding
the connection between the relatively short cell cycle durations
and the relatively long periods of peripheral cell oscillations in
some periodic hematological diseases.
\end{abstract}

\bigskip{}


\section{Introduction}

Cellular population models have been investigated intensively
since the 1960's (see, for example, Trucco \cite{trucco1965,
trucco1966}, Nooney \cite{nooney}, Rubinow \cite{rubinow1968} and
Rubinow and Lebowitz \cite{rubinow1975}) and still interest a lot
of researchers. This interest is greatly motivated, on one hand,
by the medical applications and, on the other hand, by the
biological phenomena (such as oscillations, bifurcations,
traveling waves or chaos) observed in these models and, generally
speaking, in the living world (Mackey and Glass
\cite{mackeyglass}, Mackey and Milton \cite{mackeymilton}).

Hematopoiesis is the process by which primitive stem cells
proliferate and differentiate to produce mature blood cells. It is
driven by highly coordinated patterns of gene expression under the
influence of growth factors and hormones. The regulation of
hematopoiesis is about the formation of blood cell elements in the
body. White and red blood cells and platelets are produced in the
bone marrow from where they enter the blood stream. The principal
factor stimulating red blood cell production is the hormone
produced in the kidney, called erythropoietin. About $90\%$ of the
erythropoietin is secreted by renal tubular epithelial cells when
blood is unable to deliver sufficient oxygen. When the level of
oxygen in the blood decreases this leads to a release of a
substance, which in turn causes an increase in the release of the
blood elements from the marrow. There is a feedback from the blood
to the bone marrow. Abnormalities in the feedback are considered
as major suspects in causing periodic hematological diseases, such
as auto-immune hemolytic anemia (B\'elair {\it et al.}
\cite{bmm95} and Mahaffy {\it et al.} \cite{mbm98}), cyclical
neutropenia (Haurie {\it et al.} \cite{hdm98}), chronic
myelogenous leukemia (Fowler and Mackey \cite{fm02} and
Pujo-Menjouet {\it et al.} \cite{mackeypujo2}), etc.

Cell biologists classified stem cells as proliferating cells and
resting cells (also called $G_0$-cells) (see Mackey
\cite{mackey1978, mackey1979}). Proliferating cells are committed
to undergo mitosis a certain time after their entrance into the
proliferating phase. Mackey supposed that this time of cytokinesis
is constant, that is, it is the same for all cells. Most of
committed stem cells are in the proliferating phase. The
$G_0$-phase, whose existence is known due to the works of Burns
and Tannock \cite{burnstannock}, is a quiescent stage in the
cellular development. However, it is usually believed that 95\% of
pluripotent stem cells are in the resting phase. Resting cells can
exit randomly to either entry into the proliferating phase or be
irremediably lost. Proliferating cells can also be lost by
apoptosis (programmed cell death).

The model of Mackey \cite{mackey1978} has been numerically studied
by Mackey and Rey \cite{mackey1992} and Crabb \emph{et al.}
\cite{mackey1996_1}. Computer simulations showed that there exist
strange behaviors of the stem cell population, such as
oscillations and bifurcations. Recently, Pujo-Menjouet and Mackey
\cite{pujomackey} proved the existence of a Hopf bifurcation which
causes periodic chronic myelogenous leukemia and showed the great
dependence of the model on the parameters.

In this paper, based on the model of Mackey \cite{mackey1978}, we
propose a more general model of hematopoiesis. We take into
account the fact that a cell cycle has two phases, that is, stem
cells in process are either in a resting phase or actively
proliferating. However, we do not suppose that all cells divide at
the same age, because this hypothesis is not biologically
reasonable. For example, it is believed that pluripotent stem
cells divide faster than committed stem cells, which are more
mature cells. There are strong evidences (see Bradford \emph{et
al.} \cite{bradford}) that indicate that the age of cytokinesis
$\tau$ is distributed on an interval $[\underline{\tau},
\overline{\tau}]$ with $\underline{\tau}\geq0$. Hence, we shall
assume that $\tau$ is distributed with a density $f$ supported on
an interval $[\underline{\tau}, \overline{\tau}]$, with
$0\leq\underline{\tau}<\overline{\tau}<+\infty$. The resulting
model is a system of two differential equations with distributed
delay. A simpler model, dealing with the pluripotent stem cell
population behavior, has been studied by Adimy {\it et al.}
\cite{acr04}.

Some results about stability of differential equations with
distributed delay can be mentioned. In \cite{b1989}, Boese studied
the stability of a differential equation with gamma-distributed
delay. Gamma distributions have the property to simplify the
nature of the delay and this situation is close to the one with
discrete delay. Anderson \cite{a1991,a1992} showed stability
results linked to the different moments (especially the
expectation and the variance) of the distribution. Kuang
\cite{k1994} also obtained general stability results for systems
of delay differential equations. More recently, sufficient
conditions for the stability of delay differential equations with
distributed delay have been obtained by Bernard {\it et al.}
\cite{bbm2001}. They used some properties of the distribution to
prove these results. However, in all these works, the authors
focused on sufficient conditions for the stability, there is no
necessary condition in these studies, and these results are not
applicable directly to the model considered in this paper.

This paper is organized as follows: in Section
\ref{sectionpresentation}, we present the model and establish
boundedness properties of the solutions. In Section
\ref{sectionstability}, we study the asymptotic stability of the
equilibria. We give conditions for the trivial equilibrium to be
globally asymptotically stable in Section
\ref{sectiontrivialequilibrium} and investigate the stability of
the nontrivial equilibrium in Section
\ref{sectionotherequilibrium}. In Section
\ref{sectionhopfbifurcation}, we show that a local Hopf
bifurcation occurs in our model. In Section
\ref{sectiondiscussion}, numerical simulations are performed to
demonstrate that our results can be used to explain the long
period oscillations observed in chronic myelogenous leukemia.

\section{The hematopoiesis process: presentation of the model} \label{sectionpresentation}

Denote by $r(t,a)$ and $p(t,a)$ the population densities of
resting an proliferating cells, respectively, which have spent a
time $a\geq0$ in their phase at time $t\geq0$. Resting cells can
either be lost randomly at a rate $\delta\geq0$, which takes into
account the cellular differentiation, or entry into the
proliferating phase at a rate $\beta$. Proliferating cells can be
lost by apoptosis (a programmed cell death) at a rate $\gamma \geq
0$ and, at mitosis, cells with age $a$ divide in two daughter
cells (which immediately enter the $G_0$-phase) with a rate
$g(a)$.

The function $g:[0,\overline{\tau})\to\mathbb{R}^+$ satisfies
$g(a)=0$ if $a<\underline{\tau}$ with
$0\leq\underline{\tau}<\overline{\tau}<+\infty.$ Moreover, it is
assumed to be piecewisely continuous such that
$\int_{\underline{\tau}}^{\overline{\tau}}g(a)da=+\infty$. The
later assumption describes the fact that cells which did not die
have to divide before they reach the maximal age
$\overline{\tau}$.

The nature of the trigger signal for introduction in the
proliferating phase is not clear. However, the work of Sachs
\cite{sachs} shows that we can reasonably think that it strongly
depends on the entire resting cell population, that is
$\beta=\beta(x(t))$, with
\begin{displaymath}
x(t)=\int_0^{+\infty}r(t,a)da, \qquad t\geq0.
\end{displaymath}
The function $\beta$ is supposed to be continuous and positive.
Furthermore, from a reasonable biological point of view, we assume
that $\beta$ is decreasing with $\lim_{x\to +\infty} \beta(x)=0$.
This describes the fact that the rate of re-entry into the
proliferating compartment is a decreasing function of the
$G_0$-phase population.

Usually, it is believed that the function $\beta$ is a monotone
decreasing Hill function (see Mackey \cite{mackey1978}), given by
\begin{equation}\label{beta}
\beta(x)=\beta_0\frac{\theta^n}{\theta^n+x^n}, \qquad x\geq 0,
\end{equation}
with $\beta_0>0$, $\theta\geq0$ and $n>0$. $\beta_0$ is the
maximal rate of re-entry in the proliferating phase, $\theta$ is
the number of resting cells at which $\beta$ has its maximum rate
of change with respect to the resting phase population, and $n$
describes the sensitivity of the reintroduction rate with changes
in the population.

The above parameters values are usually chosen (see Mackey
\cite{mackey1978}) to be
\begin{equation}\label{parameters}
\delta=0.05 \textrm{ day}^{-1}, \quad \gamma=0.2 \textrm{
day}^{-1}, \quad \beta_0=1.77 \textrm{ day}^{-1} \quad \textrm{
and } \quad n=3.
\end{equation}
Although an usual value of $\theta$ is $\theta=1.62 \times 10^8$
cells/kg, it can be normalized without loss of generality when one
makes a qualitative analysis of the population.

Then $r(t,a)$ and $p(t,a)$ satisfy the system of partial
differential equations
\begin{eqnarray}
\frac{\partial r}{\partial t}+\frac{\partial r}{\partial a} &=&
-\big(\delta+\beta(x(t))\big)r,\qquad a>0, \ t>0,\label{n}\vspace{1ex}\\
\frac{\partial p}{\partial t}+\frac{\partial p}{\partial a} &=&
-\big(\gamma+g(a)\big) p,\qquad\quad\ 0<a<\overline{\tau}, \ t>0,
\label{p}
\end{eqnarray}
with
\begin{displaymath}
r(0,a)=\nu(a), \ a\geq0, \qquad p(0,a)=\Gamma(a), \
a\in[0,\overline{\tau}].
\end{displaymath}
The functions $\nu=\nu(a)$ and $\Gamma=\Gamma(a)$ give the
population densities of cells which have spent a time $a$ in the
resting and proliferating phase, respectively, at time $t=0$; that
is the initial populations of cells with  age $a$ in each phase.

The boundary conditions of system (\ref{n})--(\ref{p}), which
describe the cellular flux between the two phases, are given by
\begin{displaymath}
\left\{ \begin{array}{rcl}
r(t,0)&=&2\displaystyle\int_{\underline{\tau}}^{\overline{\tau}}g(\tau)p(t,\tau)d\tau,\vspace{1ex}\\
p(t,0)&=&\beta(x(t))x(t).
\end{array}\right.
\end{displaymath}
Moreover, we suppose that $\lim_{a\to+\infty}r(t,a)=0$ and
$\lim_{a\to\overline{\tau}}p(t,a)=0$.

Let $y(t)$ denote the total population density of proliferating
cells at time $t$; then
\begin{displaymath}
y(t)=\int_0^{\overline{\tau}}p(t,a)da, \qquad t\geq0.
\end{displaymath}
Thus, integrating (\ref{n}) and (\ref{p}) with respect to the age
variable, we obtain
\begin{eqnarray}
\displaystyle\frac{dx}{dt} &=&-\big(\delta+\beta(x(t))\big)x(t)
+2\displaystyle\int_{\underline{\tau}}^{\overline{\tau}}g(\tau)p(t,\tau)d\tau,\vspace{1ex}\label{eqx}\\
\displaystyle\frac{dy}{dt} &=&-\gamma y(t)+\beta(x(t))x(t)
-\displaystyle\int_{\underline{\tau}}^{\overline{\tau}}g(\tau)p(t,\tau)d\tau.\label{eqy}
\end{eqnarray}
We define a function $G$ by
\begin{displaymath}
G(t,a)=\left\{\begin{array}{ll}
g(a)\exp\left(-\displaystyle\int_{a-t}^ag(s)ds\right),&\quad\textrm{
if }
t<a,\vspace{1ex}\\
g(a)\exp\left(-\displaystyle\int_{0}^ag(s)ds\right),&\quad\textrm{
if } a<t.
\end{array}\right.
\end{displaymath}
Set
\begin{displaymath}
f(\tau):=g(\tau)\exp\left(-\displaystyle\int_{0}^{\tau}g(s)ds\right),
\quad \tau>0.
\end{displaymath}
One can check that $f$ is a density function, supported on
$[\underline{\tau},\overline{\tau}]$, and $f$ represents the
density of division of proliferating cells. In particularly,
$\int_{\underline{\tau}}^{\overline{\tau}}f(\tau)d\tau=1$.

Using the method of characteristics to determine $p(t,a)$, we
deduce, from (\ref{eqx})--(\ref{eqy}), that the process of
hematopoiesis is described by the following system:
\begin{eqnarray} \label{restingphase}
\left\{
\begin{array}{rcl}
\displaystyle\frac{dx}{dt}&=&-\big(\delta+\beta(x(t))\big)x(t)\\
&&+\left\{
\begin{array}{ll}
2e^{-\gamma
t}\displaystyle\int_{\underline{\tau}}^{\overline{\tau}}
G(t,\tau)\Gamma(\tau-t)d\tau,& 0 \leq t \leq \underline{\tau},\vspace{1ex}\\
2\displaystyle\int_{\underline{\tau}}^t
e^{-\gamma\tau}f(\tau)\beta(x(t-\tau))
x(t-\tau) d\tau & \\
\qquad\qquad + 2e^{-\gamma t}\displaystyle\int_t^{\overline{\tau}}
G(t,\tau)\Gamma(\tau-t)d\tau,& \underline{\tau}\leq t \leq \overline{\tau},\vspace{1ex}\\
2\displaystyle\int_{\underline{\tau}}^{\overline{\tau}}
e^{-\gamma\tau}f(\tau)\beta(x(t-\tau)) x(t-\tau) d\tau,&
\overline{\tau}\leq t,
\end{array}\right. \vspace{1ex}\\
\displaystyle\frac{dy}{dt} &=&-\gamma y(t)+\beta(x(t))x(t) \vspace{1ex} \\
&& - \left\{
\begin{array}{ll}
e^{-\gamma t}
\displaystyle\int_{\underline{\tau}}^{\overline{\tau}}
G(t,\tau)\Gamma(\tau-t)d\tau,& 0 \leq t \leq \underline{\tau},\vspace{1ex}\\
\displaystyle\int_{\underline{\tau}}^t
e^{-\gamma\tau}f(\tau)\beta(x(t-\tau))x(t-\tau) d\tau & \\
\qquad\qquad + e^{-\gamma t}\displaystyle\int_t^{\overline{\tau}}
G(t,\tau)\Gamma(\tau-t)d\tau,& \underline{\tau}\leq t \leq \overline{\tau},\vspace{1ex}\\
\displaystyle\int_{\underline{\tau}}^{\overline{\tau}}
e^{-\gamma\tau}f(\tau)\beta(x(t-\tau)) x(t-\tau) d\tau,&
\overline{\tau}\leq t.
\end{array}\right.
\end{array}
\right.
\end{eqnarray}

One can give a direct biological explanation of system
(\ref{restingphase}).

In the equation for the resting cells $x(t)$, the first term in
the right-hand side accounts for $G_0$-cell loss due to either the
mortality and cellular differentiation ($\delta$) or introduction
in the proliferating phase ($\beta$). The second term represents a
cellular gain due to the movement of proliferating cells one
generation earlier. It requires some explanations. First, we
recall that all cells divide according to the density $f$,
supported on $[\underline{\tau}, \overline{\tau}]$. We shall call,
in the following, new proliferating cells, the resting cells
introduced in the proliferating phase at the considered time $t$.
When $t\leq \underline{\tau}$, no new proliferating cell is mature
enough to divide, because cells cannot divide before they have
spent a time $\underline{\tau}$ in the proliferating phase.
Therefore, the cellular gain can only proceed from cells initially
in the proliferating phase. When
$t\in[\underline{\tau},\overline{\tau}]$, the cellular increase is
obtained by division of new proliferating cells and by division of
the initial population. Finally, when $t\geq\overline{\tau}$, all
initial proliferating cells have divided or died, and the cellular
gain is obtained by division of new proliferating cells introduced
one generation earlier. The factor $2$ always accounts for the
division of each cell into two daughter cells at mitosis. The term
$e^{-\gamma t}$, with $t\in[0,\overline{\tau}]$, describes the
attenuation of the population, in the proliferating phase, due to
apoptosis.

In the equation for the proliferating cells $y(t),$ the first term
in the right-hand side accounts for cellular loss by apoptosis and
the second term is for cellular entry from the $G_0$-phase. The
last term accounts for the flux of proliferating cells to the
resting compartment.

We set $\mu:=\int_0^{\infty}\nu(a)da$. Then, initially, the
populations in the two phases are given by
\begin{displaymath}
x(0)=\mu \qquad \textrm{ and } \qquad
y(0)=\int_{0}^{\overline{\tau}}\Gamma(a)da.
\end{displaymath}
At this point, one can make a remark. Since resting cells are
introduced in the proliferating phase with a rate $\beta$, then
$\Gamma(0)$, which represents the population of cells introduced
at time $t=0$ in the cycle, must satisfy
\begin{displaymath}
\Gamma(0)=\beta(\mu)\mu.
\end{displaymath}
Taking into account the inevitable loss of proliferating cells by
apoptosis and by division, we suppose that $\Gamma(a)$ is given by
\begin{equation} \label{Gamma}
\Gamma(a)=\left\{\begin{array}{ll} e^{-\gamma a}\beta(\mu)\mu,&
\quad \textrm{ if }
a\in[0,\underline{\tau}),\vspace{1ex}\\
e^{-\gamma
a}\exp\left(-\displaystyle\int_{\underline{\tau}}^ag(s)ds\right)\beta(\mu)\mu,&
\quad \textrm{ if } a\in[\underline{\tau}, \overline{\tau}).
\end{array}\right.
\end{equation}
This simply describes that $\Gamma$ satisfies (\ref{p}) (see Webb
\cite{webb1985}, page 8). With (\ref{Gamma}) and integrating by
parts, the initial conditions of system (\ref{restingphase})
become
\begin{equation}\label{initialconditions}
x(0)=\mu, \qquad
y(0)=\beta(\mu)\mu\int_{\underline{\tau}}^{\overline{\tau}}f(\tau)\left(\frac{1-e^{-\gamma\tau}}{\gamma}\right)d\tau.
\end{equation}
When $\gamma=0$, we have
\begin{displaymath}
y(0)=\beta(\mu)\mu\int_{\underline{\tau}}^{\overline{\tau}}\tau
f(\tau)d\tau.
\end{displaymath}

Assume that the function $x\mapsto x\beta(x)$ is Lipschitz
continuous. It is immediate to show by steps that, for all
$\mu\geq0$, the system (\ref{restingphase}) under condition
(\ref{initialconditions}) has a unique nonnegative continuous
solution $(x(t),y(t))$ defined on $[0,+\infty)$.

One can notice that problem (\ref{restingphase}) reduces to a
system of two delay differential equations, with initial
conditions solutions of a system of ordinary differential
equations. On $[0,\underline{\tau}]$, the first equation for
$x(t)$ in system (\ref{restingphase}) reduces to the ordinary
differential equation
\begin{equation} \label{phitildeproblem}
\left\{ \setlength\arraycolsep{2pt}\begin{array}{rcll}
\displaystyle\frac{d\widetilde{\varphi}}{dt} &=&
-\big(\delta+\beta(\widetilde{\varphi}(t))\big)\widetilde{\varphi}(t)
+2\beta(\mu)\mu\displaystyle\int_{\underline{\tau}}^{\overline{\tau}}
e^{-\gamma\tau}f(\tau)d\tau,& \quad 0 \leq t \leq \underline{\tau},\\
\widetilde{\varphi}(0)&=&\mu, &
\end{array} \right.
\end{equation}
and, on $[\underline{\tau},\overline{\tau}]$, the second equation
reduces to the following nonautonomous delay differential equation
\begin{equation} \label{phiproblem}
\left\{ \setlength\arraycolsep{2pt}\begin{array}{rcll}
\displaystyle\frac{d\varphi}{dt} &=&
-\big(\delta+\beta(\varphi(t))\big)\varphi(t)
+ 2\beta(\mu)\mu\displaystyle\int_t^{\overline{\tau}} e^{-\gamma\tau}f(\tau)d\tau&\\
                    & & +2\displaystyle\int_{\underline{\tau}}^t
                    e^{-\gamma\tau}f(\tau)\beta(\varphi(t-\tau)) \varphi(t-\tau) d\tau,
                    & t \in[ \underline{\tau}, \overline{\tau}],\\
\varphi(t)&=&\widetilde{\varphi}(t),& t\in[0,\underline{\tau}],
\end{array} \right.
\end{equation}
where $\widetilde{\varphi}(t)$ is the unique solution of
(\ref{phitildeproblem}) for the initial condition $\mu$.

By the same way, the solution $y(t)$ of the second equation in
(\ref{restingphase}), denoted by $\psi(t)$, is given in terms of
the unique solution $\widetilde{\varphi}(t)$ of
(\ref{phitildeproblem}), associated with $\mu$, and the unique
solution $\varphi(t)$ of (\ref{phiproblem}), for
$t\in[0,\overline{\tau}]$.

Then, the system (\ref{restingphase}) can be written as an
autonomous system of delay differential equations, for
$t\geq\overline{\tau}$,
\begin{subequations}\label{problem}
\begin{eqnarray}
&& \displaystyle\frac{dx}{dt}=-\big(\delta+\beta(x(t))\big)x(t)
+2\displaystyle\int_{\underline{\tau}}^{\overline{\tau}}
e^{-\gamma\tau}f(\tau)\beta(x(t-\tau)) x(t-\tau) d\tau, \label{equationx}\\[3mm]
&&\displaystyle\frac{dy}{dt}=-\gamma y(t)+\beta(x(t))x(t)
-\displaystyle\int_{\underline{\tau}}^{\overline{\tau}}
e^{-\gamma\tau}f(\tau)\beta(x(t-\tau)) x(t-\tau)
d\tau,\label{equationy}
\end{eqnarray}
\end{subequations}
with, for $t\in[0, \overline{\tau}]$,
\begin{equation}\label{icproblem}
x(t) = \varphi(t), \quad y(t) = \psi(t).
\end{equation}
The solutions of (\ref{equationy}) are given explicitly by
\begin{equation}\label{y}
y(t)=\int_{\underline{\tau}}^{\overline{\tau}}
f(\tau)\bigg(\int_{t-\tau}^t e^{-\gamma(t-s)}\beta(x(s))x(s) \ ds
\bigg)d\tau \qquad \textrm{ for } t\geq\overline{\tau}.
\end{equation}
One can notice that $y(t)$ does not depend anymore on the initial
population $\Gamma(a)$ after one generation, that is when
$t\geq\overline{\tau}$. This can be explained as follows: Cells
initially in the proliferating phase have divided or died after
one generation; hence, new cells in the proliferating phase can
only come from resting cells $x(t)$.

On the other hand, one may have already noticed that the solutions
of (\ref{equationx}) do not depend on the solutions of
(\ref{equationy}) whereas the converse is not true. The expression
of $y(t)$ in (\ref{y}) gives more precise information on the
influence of the behavior of $x(t)$ on the stability of the
solutions $y(t)$. These results are proved in the following lemma.

\begin{lemma}\label{lemma}
Let $(x(t),y(t))$ be a solution of (\ref{problem}). If
$\lim_{t\to+\infty}x(t)$ exists and equals $C\geq0$, then
\begin{equation}\label{limity}
\lim_{t\to +\infty} y(t)=\left\{\begin{array}{ll} \beta(C)C
\displaystyle\int_{\underline{\tau}}^{\overline{\tau}}
f(\tau)\bigg(\frac{1-e^{-\gamma\tau}}{\gamma}\bigg)d\tau,&\qquad
\textrm{if } \ \gamma>0,\vspace{1ex}\\
\beta(C)C \displaystyle\int_{\underline{\tau}}^{\overline{\tau}}
\tau f(\tau)d\tau,&\qquad\textrm{if } \ \gamma=0.
\end{array}\right.
\end{equation}
If $x(t)$ is $P$-periodic, then $y(t)$ is also $P$-periodic.
\end{lemma}

\begin{proof}
By using (\ref{y}), we obtain that
\begin{equation}\label{yexplicit2}
y(t)=\int_{\underline{\tau}}^{\overline{\tau}}
f(\tau)\bigg(\int_{0}^{\tau} e^{-\gamma s}\beta(x(t-s))x(t-s) \ ds
\bigg) d\tau \qquad \textrm{ for } t\geq\overline{\tau}.
\end{equation}
Hence,
\begin{displaymath}
\lim_{t\to+\infty}y(t)=\beta(C)C\int_{\underline{\tau}}^{\overline{\tau}}
f(\tau)\bigg(\int_0^{\tau} e^{-\gamma s}\ ds\bigg) d\tau,
\end{displaymath}
and (\ref{limity}) follows immediately.

When $x(t)$ is $P$-periodic, then using (\ref{yexplicit2}) it is
obvious to see that $y(t)$ is also periodic with the same period.
\qquad
\end{proof}

Lemma \ref{lemma} shows the influence of (\ref{equationx}) on the
stability of the entire system, since the stability of solutions
of (\ref{equationx}) leads to stability of the solutions of
(\ref{equationy}).

Before studying the stability of (\ref{equationx}), we prove a
boundedness result for the solutions of this equation. The proof
is based on the one given by Mackey and Rudnicki \cite{mackey1994}
for a differential equation with a discrete delay.

\begin{proposition}\label{propbounded}
Assume that $\delta>0$. Then the solutions of (\ref{equationx})
are bounded.
\end{proposition}

\begin{proof}
Assume that $\delta>0$ and
$2(\int_{\underline{\tau}}^{\overline{\tau}}e^{-\gamma\tau}f(\tau)d\tau)\beta(0)\geq\delta$.
Since $\beta$ is decreasing and $\lim_{x\to+\infty}\beta(x)=0$,
there exists a unique $x_0\geq0$ such that
\begin{displaymath}
2\bigg(\int_{\underline{\tau}}^{\overline{\tau}}
e^{-\gamma\tau}f(\tau)d\tau\bigg)\beta(x_0)=\delta
\end{displaymath}
and
\begin{equation}\label{cond}
2\bigg(\int_{\underline{\tau}}^{\overline{\tau}}
e^{-\gamma\tau}f(\tau)d\tau\bigg)\beta(x)\leq\delta \qquad
\textrm{ for } x\geq x_0.
\end{equation}
If
$2(\int_{\underline{\tau}}^{\overline{\tau}}e^{-\gamma\tau}f(\tau)d\tau)\beta(0)<\delta$,
then (\ref{cond}) holds with $x_0=0$. Set
\begin{displaymath}
x_1:=2\bigg(\int_{\underline{\tau}}^{\overline{\tau}}
e^{-\gamma\tau}f(\tau)d\tau\bigg)\frac{\beta(0)x_0}{\delta}\geq0.
\end{displaymath}
One can check that
\begin{equation}\label{propertybeta}
2\bigg(\int_{\underline{\tau}}^{\overline{\tau}}
e^{-\gamma\tau}f(\tau)d\tau\bigg)\max_{0\leq y\leq
x}\Big(\beta(y)y\Big) \leq \delta x \qquad \textrm{ for } x\geq
x_1.
\end{equation}
Indeed, let $y\in[0,x)$. If $y\leq x_0$, then
\begin{displaymath}
2\bigg(\int_{\underline{\tau}}^{\overline{\tau}}
e^{-\gamma\tau}f(\tau)d\tau\bigg)\beta(y)y \leq
2\bigg(\int_{\underline{\tau}}^{\overline{\tau}}
e^{-\gamma\tau}f(\tau)d\tau\bigg)\beta(0)x_0=\delta x_1\leq \delta
x
\end{displaymath}
and, if $y>x_0$, then
\begin{displaymath}
2\bigg(\int_{\underline{\tau}}^{\overline{\tau}}
e^{-\gamma\tau}f(\tau)d\tau\bigg)\beta(y)y \leq \delta y\leq
\delta x.
\end{displaymath}
Hence, (\ref{propertybeta}) holds.

Assume, by contradiction, that
$\limsup_{t\to+\infty}x(t)=+\infty$, where $x(t)$ is a solution of
(\ref{equationx}). Then, there exists $t_0>\overline{\tau}$ such
that
\begin{displaymath}
x(t)\leq x(t_0) \ \textrm{ for } t\in[t_0-\overline{\tau},t_0]
\qquad \textrm{ and } \qquad x(t_0)>x_1.
\end{displaymath}
With (\ref{propertybeta}), we obtain that
\begin{displaymath}
2\displaystyle\int_{\underline{\tau}}^{\overline{\tau}}
e^{-\gamma\tau}f(\tau)\beta(x(t_0-\tau)) x(t_0-\tau) d\tau \leq
\delta x(t_0).
\end{displaymath}
This yields, with (\ref{equationx}), that
\begin{displaymath}
\frac{dx}{dt}(t_0)\leq-\beta(x(t_0))\big)x(t_0)<0,
\end{displaymath}
which gives a contradiction. Hence,
$\limsup_{t\to+\infty}x(t)<+\infty$. \qquad
\end{proof}

When $\delta=0$, the solutions of (\ref{equationx}) may not be
bounded. We show, in the next proposition, that these solutions
may explode under some conditions. However, one can notice, using
(\ref{yexplicit2}), that the solutions of (\ref{equationy}) may
still be stable in this case.

\begin{proposition} \label{propunbounded}
Assume that $\delta=0$ and
\begin{equation}\label{condexplosion}
\int_{\underline{\tau}}^{\overline{\tau}}e^{-\gamma\tau}f(\tau)d\tau
>\frac{1}{2}.
\end{equation}
In addition, assume that there exists $\overline{x}\geq0$ such
that the function $x\mapsto x\beta(x)$ is decreasing for
$x\geq\overline{x}$. If $\mu\geq\overline{x}$, then the unique
solution $x(t)$ of (\ref{equationx}) satisfies
\begin{displaymath}
\lim_{t\to+\infty} x(t)=+\infty.
\end{displaymath}
\end{proposition}

\begin{proof}
One can notice that, if $\lim_{t\to+\infty}x(t)=C$ exists, then
(\ref{equationx}) leads to
\begin{displaymath}
\bigg(2\int_{\underline{\tau}}^{\overline{\tau}}
e^{-\gamma\tau}f(\tau) d\tau -1\bigg)\beta(C)C=0.
\end{displaymath}
It follows that $C=0$.

Let $\mu\geq\overline{x}$ be given. Consider the equation
\begin{equation}\label{deltazero1}
\widetilde{\varphi}^{\prime}(t)=
2\beta(\mu)\mu\displaystyle\int_{\underline{\tau}}^{\overline{\tau}}
e^{-\gamma\tau}f(\tau)d\tau-\beta(\widetilde{\varphi}(t))\widetilde{\varphi}(t)
\quad \textrm{ for } 0 \leq t \leq \underline{\tau}
\end{equation}
with $\widetilde{\varphi}(0)=\mu$. Since the function $x\mapsto
x\beta(x)$ is decreasing for $x\geq\overline{x}$, it is immediate
that every solution $\widetilde{\varphi}(t)$ of (\ref{deltazero1})
satisfies, for $t\in[0,\underline{\tau}]$,
\begin{displaymath}
\widetilde{\varphi}^{\prime}(t)\geq
\bigg(2\int_{\underline{\tau}}^{\overline{\tau}}
e^{-\gamma\tau}f(\tau) d\tau -1\bigg)\beta(\mu)\mu>0.
\end{displaymath}
Consider now the problem
\begin{equation}\label{deltazero2}
\left\{ \setlength\arraycolsep{2pt}\begin{array}{rcll}
\varphi^{\prime}(t) &=& -\beta(\varphi(t))\varphi(t) +
2\beta(\mu)\mu\displaystyle\int_t^{\overline{\tau}}
e^{-\gamma\tau}f(\tau)d\tau&\\
                    & & +2\displaystyle\int_{\underline{\tau}}^t
                    e^{-\gamma\tau}f(\tau)\beta(\varphi(t-\tau)) \varphi(t-\tau)
                    d\tau,& \quad t \in[ \underline{\tau}, \overline{\tau}],\\
\varphi(t)&=&\widetilde{\varphi}(t),& \quad
t\in[0,\underline{\tau}],
\end{array} \right.
\end{equation}
where $\widetilde{\varphi}(t)$ is the unique solution of
(\ref{deltazero1}) for the initial condition $\mu$. Then,
\begin{displaymath}
\varphi^{\prime}(\underline{\tau})\geq
\bigg(2\int_{\underline{\tau}}^{\overline{\tau}}
e^{-\gamma\tau}f(\tau) d\tau -1\bigg)\beta(\mu)\mu>0.
\end{displaymath}
So, there exists $\varepsilon>0$ such that
$\underline{\tau}+\varepsilon\leq\overline{\tau}$ and
$\varphi^{\prime}(t)>0$ for
$t\in[\underline{\tau},\underline{\tau}+\varepsilon)$. Since
$\mu\leq
\varphi(\underline{\tau})\leq\varphi(\tau)\leq\varphi(\underline{\tau}+\varepsilon)$,
for $\tau\in[\underline{\tau},\underline{\tau}+\varepsilon]$, we
have
\begin{displaymath}
\begin{array}{rcl}
\varphi^{\prime}(\underline{\tau}+\varepsilon)&\geq&
\bigg(2\displaystyle\int_{\underline{\tau}+\varepsilon}^{\overline{\tau}}
e^{-\gamma\tau}f(\tau)d\tau-1\bigg)
\beta(\varphi(\underline{\tau}+\varepsilon))\varphi(\underline{\tau}+\varepsilon)\vspace{2ex}\\
&&+2\bigg(\displaystyle\int_{\underline{\tau}}^{\underline{\tau}+\varepsilon}
e^{-\gamma\tau}f(\tau)d\tau\bigg)\beta(\varphi(\underline{\tau}+\varepsilon))
\varphi(\underline{\tau}+\varepsilon) \vspace{2ex}\\
&\geq&\bigg(2\displaystyle\int_{\underline{\tau}}^{\overline{\tau}}
e^{-\gamma\tau}f(\tau)d\tau-1\bigg)
\beta(\varphi(\underline{\tau}+\varepsilon))\varphi(\underline{\tau}+\varepsilon).
\end{array}
\end{displaymath} Condition (\ref{condexplosion}) leads to
$\varphi^{\prime}(\underline{\tau}+\varepsilon)>0$. Using a
similar argument, we obtain that
\begin{displaymath}
\varphi^{\prime}(t)>0 \qquad \textrm{ for }
t\in[\underline{\tau},\overline{\tau}].
\end{displaymath}
To conclude, consider the delay differential equation
\begin{equation}\label{deltazero3}
x^{\prime}(t)=2\displaystyle\int_{\underline{\tau}}^{\overline{\tau}}
e^{-\gamma\tau}f(\tau)\beta(x(t-\tau)) x(t-\tau) d\tau
-\beta(x(t))x(t)
\end{equation}
with an initial condition given on
$[\underline{\tau},\overline{\tau}]$ by the solution $\varphi(t)$
of (\ref{deltazero2}). Using the same reasoning as in the previous
cases, we obtain that
\begin{displaymath}
x^{\prime}(\overline{\tau})>0.
\end{displaymath}
We thus deduce that
\begin{displaymath}
x^{\prime}(t)>0 \qquad \textrm{ for } t\geq 0.
\end{displaymath}
This completes the proof.
\end{proof}

The assumption on the function $x\mapsto x\beta(x)$ in Proposition
\ref{propunbounded} is satisfied for example when $\beta$ is given
by (\ref{beta}), with $n>1$. In this case, we can take
$\overline{x}~=~\theta/(n-~1)^{1/n}$.

We now turn our attention to the stability of (\ref{problem}).
Problem (\ref{problem}) has at most two equilibria. The first one,
$E_0=(0,0)$, always exists: it corresponds to the extinction of
the population. The second one describes the expected equilibrium
of the population; it is a nontrivial equilibrium $E^*=(x^*,y^*)$,
where $x^*$ is the unique solution of
\begin{equation}\label{Estarx}
\bigg( 2\displaystyle\int_{\underline{\tau}}^{\overline{\tau}}
e^{-\gamma\tau}f(\tau)d\tau -1 \bigg)\beta(x^*)=\delta
\end{equation}
and, from (\ref{restingphase}) and (\ref{initialconditions}),
\begin{equation}\label{Estary}
y^*=\left\{\begin{array}{ll}
\beta(x^*)x^*\displaystyle\int_{\underline{\tau}}^{\overline{\tau}}
f(\tau)\bigg(\frac{1-e^{-\gamma\tau}}{\gamma}\bigg)d\tau,
& \textrm{ if } \gamma>0,\vspace{1ex}\\
\delta x^*\displaystyle\int_{\underline{\tau}}^{\overline{\tau}}
\tau f(\tau)d\tau,& \textrm{ if } \gamma=0.
\end{array}\right.
\end{equation}
Since $\beta$ is a positive decreasing function and
$\lim_{x\to+\infty}\beta(x)=0$, then the equilibrium $E^*$ exists
if and only if
\begin{equation}\label{condexistenceequilibrium}
0<\delta <
\bigg(2\displaystyle\int_{\underline{\tau}}^{\overline{\tau}}
e^{-\gamma\tau}f(\tau)d\tau-1 \bigg)\beta(0).
\end{equation}

We shall study in Section \ref{sectionstability} the stability of
the two equilibria $E_0$ and $E^*$. From Lemma \ref{lemma}, we
only need to focus on the behavior of the equilibria of
(\ref{equationx}), that is $x\equiv 0$ and $x\equiv x^*$, to
obtain information on the behavior of the entire population.

\section{Asymptotic stability}\label{sectionstability}

We first show that $E_0$ is globally asymptotically stable when it
is the only equilibrium, and that it becomes unstable when the
nontrivial equilibrium $E^*$ appears: a transcritical bifurcation
occurs then. In a second part, we determine conditions for the
nontrivial equilibrium $E^*$ to be asymptotically stable.

\subsection{Stability of the trivial equilibrium}\label{sectiontrivialequilibrium}

In the next theorem, we give a necessary and sufficient condition
for the trivial equilibrium of (\ref{equationx}) to be globally
asymptotically stable using a Lyapunov functional. Concerning
definition and interest of Lyapunov functionals for delay
differential equations, we refer to Hale \cite{hale}.

\begin{theorem} \label{globalstabilityx}
The trivial equilibrium of the system (\ref{problem}) is globally
asymptotically stable if
\begin{equation} \label{condstabglobtrivial}
\bigg(2\displaystyle\int_{\underline{\tau}}^{\overline{\tau}}
e^{-\gamma\tau}f(\tau)d\tau-1 \bigg)\beta(0)<\delta
\end{equation}
and unstable if
\begin{equation} \label{condstabilitytrivial}
\delta<\bigg(2\displaystyle\int_{\underline{\tau}}^{\overline{\tau}}
e^{-\gamma\tau}f(\tau)d\tau-1 \bigg)\beta(0).
\end{equation}
\end{theorem}

\begin{proof}
We first assume that (\ref{condstabglobtrivial}) holds. Denote by
$C^+$ the set of continuous nonnegative functions on
$[0,\overline{\tau}]$ and define the mapping $J : C^+ \to
[0,+\infty)$ by
\begin{displaymath}
J(\varphi)=B(\varphi(\overline{\tau}))+\int_{\underline{\tau}}^{\overline{\tau}}
e^{-\gamma\tau}f(\tau)\bigg(
\int_{\overline{\tau}-\tau}^{\overline{\tau}}\Big(\beta\big(\varphi(\theta)\big)\varphi(\theta)
\Big)^2 d\theta \bigg) d\tau
\end{displaymath}
for all $\varphi\in C^+$, where
\begin{displaymath}
B(x)=\displaystyle\int_{0}^{x} \beta(s)s\ ds \qquad \textrm{ for
all } x\geq0.
\end{displaymath}
We set (see \cite{hale})
\begin{displaymath}
\dot
J(\varphi)=\limsup_{t\to0^+}\frac{J(x^{\varphi}_t)-J(\varphi)}{t}\qquad
\textrm{ for } \varphi\in C^+,
\end{displaymath}
where $x^{\varphi}$ is the unique solution of (\ref{equationx})
associated with the initial condition $\varphi\in C^+$, and
$x^{\varphi}_t(\theta)=x^{\varphi}(t+\theta)$ for
$\theta\in[0,\overline{\tau}]$. Then,
\begin{equation}\label{derivativeJ}
\setlength\arraycolsep{2pt}
\begin{array}{rcl} \dot J(\varphi) &=&
\displaystyle\frac{d\varphi}{dt}(\overline{\tau})\beta\big(\varphi(\overline{\tau})\big)
\varphi(\overline{\tau})\vspace{1ex}\\
&&+\displaystyle\int_{\underline{\tau}}^{\overline{\tau}}
e^{-\gamma\tau}f(\tau)\bigg(\Big(\beta\big(\varphi(\overline{\tau})\big)\varphi(\overline{\tau})
\Big)^2-\Big(\beta\big(\varphi(\overline{\tau}-\tau)\big)\varphi(\overline{\tau}-\tau)
\Big)^2 \bigg)d\tau.
\end{array}
\end{equation}
Using (\ref{equationx}), we have
\begin{displaymath}
\displaystyle\frac{d\varphi}{dt}(\overline{\tau}) =
-\Big(\delta+\beta\big(\varphi(\overline{\tau})\big)\Big)\varphi(\overline{\tau})
+ 2\int_{\underline{\tau}}^{\overline{\tau}}
e^{-\gamma\tau}f(\tau)\beta\big(\varphi(\overline{\tau}-\tau)\big)\varphi(\overline{\tau}-\tau)d\tau.
\end{displaymath}
Therefore, (\ref{derivativeJ}) becomes
\begin{eqnarray}
\dot J(\varphi) & = &
-\Big(\delta+\beta\big(\varphi(\overline{\tau})\big)\Big)
\beta\big(\varphi(\overline{\tau})\big)\varphi^{2}(\overline{\tau})
+\int_{\underline{\tau}}^{\overline{\tau}}e^{-\gamma\tau}f(\tau)
\bigg[\Big(\beta\big(\varphi(\overline{\tau})\big)\varphi(\overline{\tau}) \Big)^2 \nonumber \\
        &   & +2\beta\big(\varphi(\overline{\tau})\big)\varphi(\overline{\tau})
        \beta\big(\varphi(\overline{\tau}-\tau)\big)\varphi(\overline{\tau}-\tau)
        -\Big(\beta\big(\varphi(\overline{\tau}-\tau)\big)\varphi(\overline{\tau}-\tau) \Big)^2 \bigg]d\tau
        \nonumber \\
        & = & -\Big(\delta+\beta\big(\varphi(\overline{\tau})\big)\Big)
        \beta\big(\varphi(\overline{\tau})\big)\varphi^{2}(\overline{\tau})
        + 2\Big(\beta\big(\varphi(\overline{\tau})\big)\varphi(\overline{\tau}) \Big)^2
        \int_{\underline{\tau}}^{\overline{\tau}} e^{-\gamma\tau}f(\tau) d\tau \nonumber \\
        &   & -\int_{\underline{\tau}}^{\overline{\tau}}e^{-\gamma\tau}f(\tau)
        \Big[\beta\big(\varphi(\overline{\tau})\big)\varphi(\overline{\tau})
        -\beta\big(\varphi(\overline{\tau}-\tau)\big)\varphi(\overline{\tau}-\tau) \Big]^2 d\tau.
        \nonumber
\end{eqnarray}
Hence,
\begin{displaymath}
\dot J(\varphi) \leq  -u(\varphi(\overline{\tau})),
\end{displaymath}
where the function $u$ is defined, for $x \geq 0$, by
\begin{equation}\label{functionu}
u(x)=r(x)\beta(x)x^{2}
\end{equation}
with
\begin{displaymath}
r(x)=\delta-\bigg(2\displaystyle\int_{\underline{\tau}}^{\overline{\tau}}
e^{-\gamma\tau}f(\tau) d\tau-1\bigg)\beta(x).
\end{displaymath}
Since $\beta$ is decreasing, $r$ is a monotone function. Moreover,
(\ref{condstabglobtrivial}) leads to $r(0)>0$, and $\lim_{x \to
\infty} r(x)= \delta\geq 0$. Therefore, $r$ is positive on
$[0,+\infty)$.

Consequently, the function $u$ defined by (\ref{functionu}) is
nonnegative on $[0,+\infty)$ and $u(x)=0$ if and only if $x=0$. We
deduce that every solution of (\ref{equationx}), with $\varphi\in
C^{+}$, tends to zero as $t$ tends to $+\infty$.

We suppose now that (\ref{condstabilitytrivial}) holds. The
linearization of (\ref{equationx}) around $x\equiv 0$ leads to the
characteristic equation
\begin{equation}\label{characteristicequationzero}
\Delta_0(\lambda):=\lambda+\delta+\beta(0)-2\beta(0)
\displaystyle\int_{\underline{\tau}}^{\overline{\tau}}
e^{-(\lambda+\gamma)\tau}f(\tau) d\tau=0.
\end{equation}
We consider $\Delta_0$ as a real function. Since
\begin{displaymath}
\frac{d\Delta_0}{d\lambda}=1+2\beta(0)\displaystyle\int_{\underline{\tau}}^{\overline{\tau}}
\tau e^{-(\lambda+\gamma)\tau}f(\tau) d\tau >0,
\end{displaymath}
it follows that $\Delta_0$ is an increasing function. Moreover,
(\ref{characteristicequationzero}) yields
\begin{displaymath}
\lim_{\lambda\to -\infty}\Delta_0(\lambda)=-\infty, \qquad
\lim_{\lambda\to +\infty}\Delta_0(\lambda)=+\infty,
\end{displaymath}
and (\ref{condstabilitytrivial}) implies that
\begin{displaymath}
\Delta_0(0)=\delta-\bigg(2\displaystyle\int_{\underline{\tau}}^{\overline{\tau}}
e^{-\gamma\tau}f(\tau) d\tau-1 \bigg)\beta(0)< 0.
\end{displaymath}
Hence, $\Delta_0(\lambda)$ has a unique real root which is
positive. Consequently, (\ref{characteristicequationzero}) has at
least one characteristic root with positive real part. Therefore,
the equilibrium $x\equiv 0$ of (\ref{equationx}) is not stable.
This completes the proof. \qquad
\end{proof}

The inequality (\ref{condstabglobtrivial}) is satisfied when
$\delta$ or $\gamma$ (the mortality rates) is large or when
$\beta(0)$ is small. Biologically, these conditions correspond to
a population which cannot survive, because the mortality rates are
too large or, simply, because not enough cells are introduced in
the proliferating phase and, then, the population renewal is not
supplied.

\emph{Remark 1.} One can notice that when
\begin{displaymath}
\int_{\underline{\tau}}^{\overline{\tau}}e^{-\gamma\tau}f(\tau)d\tau<\frac{1}{2},
\end{displaymath}
the trivial equilibrium $E_0$ is the only equilibrium of
(\ref{problem}) and is globally asymptotically stable. When
\begin{displaymath}
\int_{\underline{\tau}}^{\overline{\tau}}e^{-\gamma\tau}f(\tau)d\tau=\frac{1}{2},
\end{displaymath}
then $E_0$ is globally asymptotically stable if $\delta>0$. When
the equality
\begin{displaymath}
\bigg(2\displaystyle\int_{\underline{\tau}}^{\overline{\tau}}
e^{-\gamma\tau}f(\tau)d\tau-1 \bigg)\beta(0)=\delta
\end{displaymath}
holds, one can check that $\lambda=0$ is a characteristic root of
(\ref{characteristicequationzero}) and all other characteristic
roots have negative real parts. Hence, we cannot conclude on the
stability or instability of the trivial equilibrium $E_0$ of
(\ref{problem}) without further analysis. However, this is not the
subject of this paper.

\subsection{Stability of the nontrivial equilibrium} \label{sectionotherequilibrium}

We concentrate, in this section, on the equilibrium
$E^*=(x^*,y^*)$ defined by (\ref{Estarx})--(\ref{Estary}). Hence,
throughout this section, we assume that
(\ref{condexistenceequilibrium}) holds, that is
\begin{displaymath}
0<\delta <
\bigg(2\displaystyle\int_{\underline{\tau}}^{\overline{\tau}}
e^{-\gamma\tau}f(\tau)d\tau-1 \bigg)\beta(0).
\end{displaymath}
Since $\delta>0$ and $\beta(0)>0$,
(\ref{condexistenceequilibrium}) implies, in particularly, that
\begin{equation}\label{condpositivityint}
\displaystyle\int_{\underline{\tau}}^{\overline{\tau}}
e^{-\gamma\tau}f(\tau)d\tau>\frac{1}{2}.
\end{equation}
From Lemma \ref{lemma}, we only need to focus on the stability of
the nontrivial equilibrium $x\equiv x^*$ of (\ref{equationx}). To
that aim, we linearize (\ref{equationx}) around $x^*$. Denote by
$\beta^*\in\mathbb{R}$ the quantity
\begin{equation}\label{betastar}
\beta^*:=\frac{d}{dx}\Big(x\beta(x)\Big)\Big|_{x=x^*}=\beta(x^*)+x^*\beta^{\prime}(x^*)
\end{equation}
and set $u(t)=x(t)-x^*$. The linearization of (\ref{equationx}) is
given by
\begin{displaymath}
\displaystyle\frac{du}{dt}=
-(\delta+\beta^*)u(t)+2\beta^*\displaystyle\int_{\underline{\tau}}^{\overline{\tau}}
e^{-\gamma\tau}f(\tau)u(t-\tau) d\tau.
\end{displaymath}
Then, the characteristic equation is
\begin{equation}\label{characteristicequation}
\Delta(\lambda):=\lambda+\delta+\beta^*-2\beta^*\displaystyle\int_{\underline{\tau}}^{\overline{\tau}}
e^{-(\lambda+\gamma)\tau}f(\tau) d\tau =0.
\end{equation}

One can notice that the function $x\mapsto x\beta(x)$ is usually
not monotone. For example, if $\beta$ is given by (\ref{beta})
with $n>1$, the function $x\mapsto x\beta(x)$ is increasing for
$x\leq \theta/(n-1)^{1/n}$ and decreasing for
$x>\theta/(n-1)^{1/n}$. In this case, $\beta^*$ is nonnegative
when $x^*$ is close to zero and negative when $x^*$ is large
enough.

The following theorem deals with the asymptotic stability of
$E^*$.

\begin{theorem} \label{theoremstability}
Assume that (\ref{condexistenceequilibrium}) holds. If
\begin{equation}\label{secondconditionbeta}
\beta^*\geq-\frac{\delta}{2\displaystyle\int_{\underline{\tau}}^{\overline{\tau}}
e^{-\gamma\tau}f(\tau)d\tau+1},
\end{equation}
then $E^*$ is locally asymptotically stable.
\end{theorem}

\begin{proof}
We first prove that the equilibrium $x\equiv x^*$ is locally
asymptotically stable when $\beta^*\geq0$. We consider the mapping
$\Delta(\lambda)$, given by (\ref{characteristicequation}), as a
real function of $\lambda$. Then $\Delta(\lambda)$ is continuously
differentiable on $\mathbb{R}$ and its first derivative is given
by
\begin{equation}\label{firstderivative}
\frac{d\Delta}{d\lambda}=1+2\beta^{*}\displaystyle\int_{\underline{\tau}}^{\overline{\tau}}
\tau e^{-(\lambda+\gamma)\tau}f(\tau) d\tau >0.
\end{equation}
Hence, $\Delta(\lambda)$ is an increasing function of $\lambda$
satisfying
\begin{displaymath}
\lim_{\lambda\to -\infty}\Delta(\lambda)=-\infty \qquad \textrm{
and } \qquad \lim_{\lambda\to +\infty} \Delta(\lambda)=+\infty.
\end{displaymath}
Then, there exists a unique $\lambda_{0}\in\mathbb{R}$ such that
$\Delta(\lambda_{0})=0$. Moreover, since
\begin{displaymath}
\Delta(0)=\delta-\bigg(2\displaystyle\int_{\underline{\tau}}^{\overline{\tau}}
e^{-\gamma\tau}f(\tau) d\tau -1\bigg)\beta^*,
\end{displaymath}
we deduce, by using (\ref{Estarx}), (\ref{condpositivityint}) and
(\ref{betastar}), that
\begin{displaymath}
\Delta(0)=-\bigg(2\int_{\underline{\tau}}^{\overline{\tau}}
e^{-\gamma\tau}f(\tau)d\tau-1\bigg)x^*\beta^{\prime}(x^*)>0.
\end{displaymath}
Consequently, $\lambda_0<0$.

Let $\lambda=\mu+i\omega$ be a characteristic root of
(\ref{characteristicequation}) such that $\mu>\lambda_{0}$.
Considering the real part of (\ref{characteristicequation}), we
obtain that
\begin{equation} \label{systemecomplet}
\mu=-(\delta+\beta^*)+2\beta^*\displaystyle\displaystyle\int_{\underline{\tau}}^{\overline{\tau}}
e^{-(\mu+\gamma)\tau}f(\tau)\cos(\omega\tau)d\tau.
\end{equation}
Using (\ref{characteristicequation}), with $\lambda=\lambda_0$,
together with (\ref{systemecomplet}), we then obtain
\begin{displaymath}
\mu-\lambda_{0}=2\beta^{*}\displaystyle\int_{\underline{\tau}}^{\overline{\tau}}
e^{-\gamma\tau}f(\tau)\big[e^{-\mu\tau}\cos(\omega\tau)-e^{-\lambda_{0}\tau}\big]d\tau.
\end{displaymath}
However,
\begin{displaymath}
e^{-\mu\tau}\cos(\omega\tau)-e^{-\lambda_{0}\tau}<0
\end{displaymath}
for all $\tau\in[\underline{\tau},\overline{\tau}]$. So we obtain
that $\mu-\lambda_{0}<0$, which leads to a contradiction. This
implies that all characteristic roots of
(\ref{characteristicequation}) have negative real part and the
equilibrium $x\equiv x^*$ of (\ref{equationx}) is locally
asymptotically stable.

Now, assume that $\beta^*<0$ and
\begin{equation}\label{secondconditionbetabis}
\beta^*>-\frac{\delta}{2\displaystyle\int_{\underline{\tau}}^{\overline{\tau}}
e^{-\gamma\tau}f(\tau)d\tau+1}.
\end{equation}
Let $\lambda=\mu+i\omega$ be a characteristic root of
(\ref{characteristicequation}) such that $\mu>0$. Since
\begin{displaymath}
\int_{\underline{\tau}}^{\overline{\tau}}
e^{-\gamma\tau}f(\tau)\Big(e^{-\mu\tau}\cos(\omega\tau)+1\Big)d\tau\geq0,
\end{displaymath}
we have
\begin{displaymath}
2\beta^*\int_{\underline{\tau}}^{\overline{\tau}}
e^{-(\mu+\gamma)\tau}f(\tau)\cos(\omega\tau)d\tau\leq-2\beta^*\int_{\underline{\tau}}^{\overline{\tau}}
e^{-\gamma\tau}f(\tau)d\tau.
\end{displaymath}
So, (\ref{systemecomplet}) and (\ref{secondconditionbetabis}) lead
to
\begin{displaymath}
\mu\leq-(\delta+\beta^*)-2\beta^*\int_{\underline{\tau}}^{\overline{\tau}}
e^{-\gamma\tau}f(\tau)d\tau<0,
\end{displaymath}
a contradiction. Therefore, $\mu\leq 0$.

Suppose now that (\ref{characteristicequation}) has a purely
imaginary characteristic root $i\omega$, with
$\omega\in\mathbb{R}$. Then, (\ref{systemecomplet}) leads to
\begin{displaymath}
\int_{\underline{\tau}}^{\overline{\tau}}
e^{-\gamma\tau}f(\tau)\cos(\omega\tau)d\tau =
\frac{\delta+\beta^*}{2\beta^*}.
\end{displaymath}
However,
\begin{displaymath}
\bigg|\int_{\underline{\tau}}^{\overline{\tau}}
e^{-\gamma\tau}f(\tau)\cos(\omega\tau)d\tau
\bigg|\leq\int_{\underline{\tau}}^{\overline{\tau}}
e^{-\gamma\tau}f(\tau)d\tau
\end{displaymath}
and (\ref{secondconditionbetabis}) yields
\begin{displaymath}
\frac{\delta+\beta^*}{2\beta^*}<-\int_{\underline{\tau}}^{\overline{\tau}}
e^{-\gamma\tau}f(\tau)d\tau.
\end{displaymath}
Hence, (\ref{characteristicequation}) has no purely imaginary
root. Consequently, all characteristic roots of
(\ref{characteristicequation}) have negative real part and the
nontrivial equilibrium $x\equiv x^*$ of (\ref{equationx}) is
locally asymptotically stable.

Finally, assume that
\begin{equation}\label{secondconditionbetater}
\beta^*=-\frac{\delta}{2\displaystyle\int_{\underline{\tau}}^{\overline{\tau}}
e^{-\gamma\tau}f(\tau)d\tau+1}.
\end{equation}
Consider a characteristic root $\lambda=\mu+i\omega$ of
(\ref{characteristicequation}), which reduces, with
(\ref{secondconditionbetater}), to
\begin{equation}\label{reducedcharacteristicequation}
\lambda-2\beta^*\int_{\underline{\tau}}^{\overline{\tau}}
e^{-\gamma\tau}f(\tau)\Big(1+e^{-\lambda\tau}\Big)d\tau=0.
\end{equation}
Suppose, by contradiction, that $\mu>0$. By considering the real
part of (\ref{reducedcharacteristicequation}), we have
\begin{displaymath}
\mu=2\beta^*\int_{\underline{\tau}}^{\overline{\tau}}
e^{-\gamma\tau}f(\tau)\Big(1+e^{-\mu\tau}\cos(\omega\tau)\Big)d\tau<0.
\end{displaymath}
We obtain a contradiction, therefore $\mu\leq0$. If we suppose now
that $\mu=0$ then we easily obtain that
\begin{displaymath}
\cos(\omega\tau)=-1 \qquad \textrm{for all }
\tau\in[\underline{\tau},\overline{\tau}],
\end{displaymath}
which is impossible. It follows that all characteristic roots of
(\ref{characteristicequation}) have negative real parts when
(\ref{secondconditionbetater}) holds and the equilibrium $x\equiv
x^*$ is locally asymptotically stable.

From Lemma \ref{lemma}, we conclude that $E^*$ is locally
asymptotically stable when (\ref{secondconditionbeta}) holds.
\qquad
\end{proof}

\begin{figure}[!htp]
\begin{center}\includegraphics[width=10cm, height=8cm]{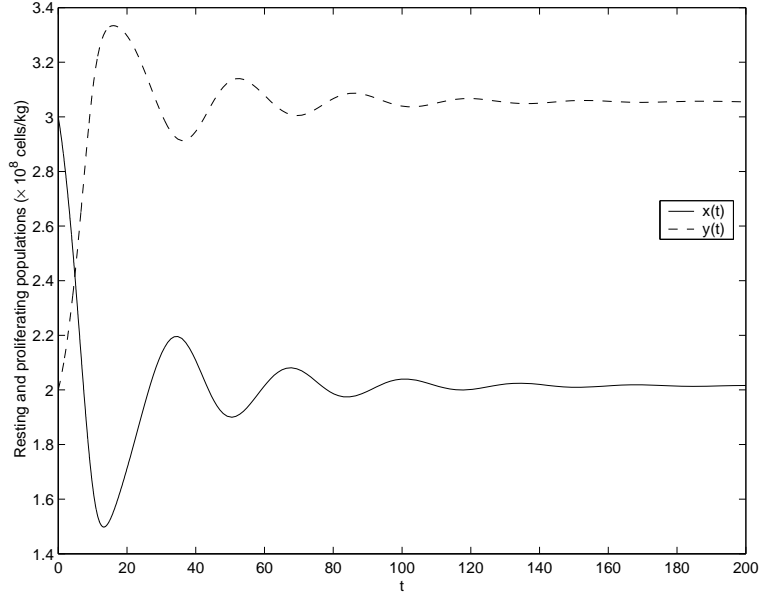}\end{center}
\caption{The solutions $x(t)$ (solid curve) and $y(t)$ (dashed
curve) of system (\ref{problem}) are drawn for values of the
parameters $\beta_0$, $\delta$ and $\gamma$ given by
(\ref{parameters}), $n=2.42$, $\underline{\tau}=0$ and
$\overline{\tau}=7$ days. In this case, the nontrivial equilibrium
$E^*$ is locally asymptotically stable, even though the solutions
oscillate transiently.}\label{stabilityfig}
\end{figure}

The asymptotic stability of $E^*$ is shown in Fig.
\ref{stabilityfig}. Values of the parameters are given by
(\ref{parameters}), except $n=2.42$, $\underline{\tau}=0$ and
$\overline{\tau}=7$ days. The function $f$ is defined by
\begin{equation}\label{uniformlaw}
f(\tau)=\left\{\begin{array}{ll}
\displaystyle\frac{1}{\overline{\tau}-\underline{\tau}},&\textrm{if
} \tau\in[\underline{\tau},\overline{\tau}],\\
0,&\textrm{otherwise.}
\end{array}\right.
\end{equation}

The {\sc Matlab} solver for delay differential equations, dde23
\cite{dde23}, is used to obtain Fig. \ref{stabilityfig}, as well
as illustrations in Section  \ref{sectionhopfbifurcation} and
Section \ref{sectiondiscussion}.

When (\ref{secondconditionbeta}) does not hold, we have
necessarily $\beta^*<0$. In this case, we cannot obtain the
stability of $E^*$ for all values of $\beta^*$. In fact, in the
next section we are going to show that the equilibrium $E^*$ can
be destabilized, in this case, via a Hopf bifurcation.

\section{Hopf bifurcation and periodic solutions} \label{sectionhopfbifurcation}

In this section, we show that the equilibrium $x\equiv x^*$ of
(\ref{equationx}) can become unstable when
(\ref{secondconditionbeta}) does not hold anymore. Throughout this
section, we assume that
\begin{displaymath}
\underline{\tau}=0
\end{displaymath}
and (\ref{condexistenceequilibrium}) holds, that is
\begin{displaymath}
0<\delta < \bigg(2\displaystyle\int_0^{\overline{\tau}}
e^{-\gamma\tau}f(\tau)d\tau-1 \bigg)\beta(0).
\end{displaymath}
From Proposition \ref{propbounded}, the solutions of
(\ref{equationx}) are bounded. Consequently, instability in
(\ref{equationx}) occurs only via oscillatory solutions.

We assume that
\begin{equation}\label{condbetastarinstability}
\beta^*<-\frac{\delta}{2\displaystyle\int_0^{\overline{\tau}}
e^{-\gamma\tau}f(\tau)d\tau+1}:=\widetilde{\delta}.
\end{equation}
Otherwise, the nontrivial equilibrium $x\equiv x^*$ of
(\ref{equationx}) is locally asymptotically stable (see Theorem
\ref{theoremstability}).

If instability occurs for a particular value
$\beta^*<\widetilde{\delta}$, a characteristic root of
(\ref{characteristicequation}) must intersect the imaginary axis.
Hence, we look for purely imaginary characteristic roots
$i\omega$, $\omega\in\mathbb{R}$, of
(\ref{characteristicequation}). If $i\omega$ is a characteristic
root of (\ref{characteristicequation}), then $\omega$ is a
solution of the system
\begin{equation}\label{systempureimaginaryroot}
\left\{ \begin{array}{rcl}
\delta+\beta^*(1-2C(\omega))&=&0,\\
\omega +2\beta^*S(\omega)&=&0,
\end{array}\right.
\end{equation}
where
\begin{displaymath}
C(\omega):=\int_0^{\overline{\tau}}
e^{-\gamma\tau}f(\tau)\cos(\omega\tau)d\tau \quad \textrm{ and }
\quad S(\omega):=\int_0^{\overline{\tau}}
e^{-\gamma\tau}f(\tau)\sin(\omega\tau)d\tau.
\end{displaymath}
One can notice that $\omega=0$ is not a solution of
(\ref{systempureimaginaryroot}). Otherwise,
\begin{displaymath}
\delta=\bigg(2\int_0^{\overline{\tau}}
e^{-\gamma\tau}f(\tau)d\tau-1\bigg)\beta^*<0,
\end{displaymath}
which gives a contradiction. Moreover, if $\omega$ is a solution
of (\ref{systempureimaginaryroot}), then $-i\omega$ is also a
characteristic root. Thus, we only look for positive solutions
$\omega$.

\begin{lemma}\label{lemmaexistenceproperties}
Assume that the function $\tau\mapsto e^{-\gamma\tau}f(\tau)$ is
decreasing. Then, for each $\delta$ such that
(\ref{condexistenceequilibrium}) is satisfied,
(\ref{systempureimaginaryroot}) has at least one solution
$(\beta^*_c,\omega_c)$, with $\beta^*_c<\widetilde{\delta}$ and
$\omega_c>0$. It follows that (\ref{characteristicequation}) has
at least one pair of purely imaginary roots $\pm i\omega_c$ for
$\beta^*=\beta^*_c$. Moreover, $\pm i\omega_c$ are simple
characteristic roots of (\ref{characteristicequation}). Consider
the branch of characteristic roots $\lambda(-\beta^*)$, such that
$\lambda(-\beta^*_c)=i\omega_c$. Then
\begin{equation}\label{conditionhopfbifurcation}
\frac{d Re(\lambda)}{d(-\beta^*)}\bigg|_{\beta^*=\beta^*_c}>0
\quad \textrm{ if and only if } \quad
-\delta\bigg(\frac{S(\omega_c)}{\omega_c}\bigg)^{\prime}>C^{\prime}(\omega_c).
\end{equation}
\end{lemma}

\begin{proof}
First, we show by induction that $S(\omega)>0$ for $\omega>0$. It
is clear that $S(\omega)>0$ if $\omega\overline{\tau}\in(0,\pi]$.
Suppose that $\omega\overline{\tau}\in(\pi,2\pi]$. Then
\begin{displaymath}
\begin{array}{rcl}
S(\omega)&=&\displaystyle\frac{1}{\omega}\int_0^{\omega\overline{\tau}}
e^{-\gamma\frac{\tau}{\omega}}f\Big(\frac{\tau}{\omega}\Big)\sin(\tau)d\tau\vspace{1ex}\\
&=&\displaystyle\frac{1}{\omega}\int_0^{\pi}e^{-\gamma\frac{\tau}{\omega}}
f\Big(\frac{\tau}{\omega}\Big)\sin(\tau)d\tau+\displaystyle\frac{1}{\omega}
\int_{\pi}^{\omega\overline{\tau}}e^{-\gamma\frac{\tau}{\omega}}
f\Big(\frac{\tau}{\omega}\Big)\sin(\tau)d\tau.
\end{array}
\end{displaymath}
Since $f$ is supported on the interval $[0,\overline{\tau}]$, it
follows that
\begin{displaymath}
\displaystyle\int_{\omega\overline{\tau}}^{2\pi}e^{-\gamma\frac{\tau}{\omega}}
f\Big(\frac{\tau}{\omega}\Big)\sin(\tau)d\tau=0.
\end{displaymath}
So, we obtain
\begin{displaymath}
\begin{array}{rcl}
S(\omega)&=&\displaystyle\frac{1}{\omega}\int_0^{\pi}e^{-\gamma\frac{\tau}{\omega}}
f\Big(\frac{\tau}{\omega}\Big)\sin(\tau)d\tau+\displaystyle\frac{1}{\omega}
\int_{\pi}^{2\pi}e^{-\gamma\frac{\tau}{\omega}}
f\Big(\frac{\tau}{\omega}\Big)\sin(\tau)d\tau\vspace{1ex}\\
&=&\displaystyle\frac{1}{\omega}\int_0^{\pi}\bigg(e^{-\gamma\frac{\tau}{\omega}}
f\Big(\frac{\tau}{\omega}\Big)-e^{-\gamma\frac{\tau+\pi}{\omega}}
f\Big(\frac{\tau+\pi}{\omega}\Big)\bigg)\sin(\tau)d\tau.
\end{array}
\end{displaymath}
Since the function $\tau\mapsto e^{-\gamma\tau}f(\tau)$ is
decreasing, we finally get $S(\omega)>0$. Using a similar argument
for $\omega\overline{\tau}\in(k\pi,(k+1)\pi]$, with
$k\in\mathbb{N}$, $k\geq2$, we deduce that $S(\omega)>0$ for all
$\omega>0$.

Consider the equation
\begin{equation}\label{equationg}
g(\omega):=\frac{\omega\big(1-2C(\omega)\big)}{2S(\omega)}=\delta,
\qquad \omega>0 .
\end{equation}
The function $g$ is continuous with
\begin{equation}
\lim_{\omega\to
0}g(\omega)=\frac{1-2C(0)}{2\displaystyle\int_0^{\overline{\tau}}
\tau e^{-\gamma\tau}f(\tau)d\tau}<0
\end{equation}
because (\ref{condexistenceequilibrium}) leads to $1-2C(0)<0$.
Moreover, the Riemann--Lebesgue's lemma implies that
\begin{displaymath}
\lim_{\omega\to+\infty}C(\omega)=\lim_{\omega\to+\infty}S(\omega)=0.
\end{displaymath}
This yields
\begin{displaymath}
\lim_{\omega\to+\infty}g(\omega)=+\infty.
\end{displaymath}
We conclude that there exists a solution $\omega_c>0$ of
(\ref{equationg}). Since $S(\omega_c)>0$ and
$g(\omega_c)=\delta>0$, we obtain $1-2C(\omega_c)>0$. Set
\begin{equation}\label{betastarcrit}
\beta^*_c=-\frac{\delta}{1-2C(\omega_c)}<0.
\end{equation}
Since $|C(\omega_c)|<C(0)$, it follows that
\begin{displaymath}
\beta^*_c<-\frac{\delta}{2C(0)+1}=\widetilde{\delta}.
\end{displaymath}
One can check that $(\beta^*_c,\omega_c)$ is a solution of
(\ref{systempureimaginaryroot}). It follows that $\pm i\omega_c$
are characteristic roots of (\ref{characteristicequation}) for
$\beta^*=\beta^*_c$.

Define a branch of characteristic roots $\lambda(-\beta^*)$ of
(\ref{characteristicequation}) such that
$\lambda(-\beta^*_c)=i\omega_c$. We use the parameter $-\beta^*$
because $\beta^*<\widetilde{\delta}<0$.

Using (\ref{characteristicequation}), we obtain
\begin{equation}\label{derivative}
\bigg[1+2\beta^*\int_0^{\overline{\tau}}\tau
e^{-(\lambda+\gamma)\tau} f(\tau)
d\tau\bigg]\frac{d\lambda}{d(-\beta^*)}
=1-2\int_0^{\overline{\tau}}e^{-(\lambda+\gamma)\tau}f(\tau)d\tau.
\end{equation}
If we assume, by contradiction, that $i\omega_c$ is not a simple
root of (\ref{characteristicequation}), then (\ref{derivative})
leads to
\begin{displaymath}
C(\omega_c)=\frac{1}{2} \qquad \textrm{ and } \qquad
S(\omega_c)=0.
\end{displaymath}
Since $S(\omega_c)>0$, we obtain a contradiction. Thus,
$i\omega_c$ is a simple root of (\ref{characteristicequation}).

Moreover, using (\ref{derivative}), we have
\begin{displaymath}
\bigg(\frac{d\lambda}{d(-\beta^*)}\bigg)^{-1}=
\frac{1+2\beta^*\displaystyle\int_0^{\overline{\tau}}\tau
e^{-(\lambda+\gamma)\tau} f(\tau)
d\tau}{1-2\displaystyle\int_0^{\overline{\tau}}e^{-(\lambda+\gamma)\tau}f(\tau)d\tau}.
\end{displaymath}
Since $\lambda$ is a characteristic root of
(\ref{characteristicequation}), we also have
\begin{displaymath}
1-2\displaystyle\int_0^{\overline{\tau}}e^{-(\lambda+\gamma)\tau}
f(\tau)d\tau=-\frac{\lambda+\delta}{\beta^*}.
\end{displaymath}
So, we deduce
\begin{displaymath}
\bigg(\frac{d\lambda}{d(-\beta^*)}\bigg)^{-1}=
-\beta^*\frac{1+2\beta^*\displaystyle\int_0^{\overline{\tau}}\tau
e^{-(\lambda+\gamma)\tau} f(\tau) d\tau}{\lambda+\delta}.
\end{displaymath}
Then,
\begin{displaymath}
\setlength\arraycolsep{2pt}
\begin{array}{rcl}
sign\bigg\{\displaystyle\frac{dRe(\lambda)}{d(-\beta^*)}\bigg\}\bigg|_{\beta^*=\beta^*_c}
&=&
sign\bigg\{Re\bigg(\displaystyle\frac{d\lambda}{d(-\beta^*)}\bigg)^{-1}\bigg\}
\bigg|_{\beta^*=\beta^*_c} \vspace{1ex}\\
&=&sign\bigg\{Re\bigg(-\beta^*\displaystyle\frac{1+2\beta^*\displaystyle\int_0^{\overline{\tau}}\tau
e^{-(\lambda+\gamma)\tau} f(\tau)
d\tau}{\lambda+\delta}\bigg)\bigg\}\bigg|_{\beta^*=\beta^*_c}
\vspace{1ex}\\
&=&sign\bigg\{-\beta^*_c\displaystyle\frac{\delta(1+2\beta^*_cS^{\prime}(\omega_c))
+2\beta^*_c\omega_cC^{\prime}(\omega_c)}{\delta^2+\omega_c^2}\bigg\} \vspace{1ex}\\
&=&sign\bigg\{\delta(1+2\beta^*_cS^{\prime}(\omega_c))+2\beta^*_c\omega_c
C^{\prime}(\omega_c)\bigg\}.
\end{array}
\end{displaymath}
From (\ref{betastarcrit}) and the fact that $1-2C(\omega_c)>0$,
this leads to
\begin{displaymath}
\begin{array}{rcl}
sign\bigg\{\displaystyle\frac{dRe(\lambda)}{d(-\beta^*)}\bigg\}\bigg|_{\beta^*=\beta^*_c}&=&
sign\bigg\{1-2C(\omega_c)-2\delta
S^{\prime}(\omega_c)-2\omega_cC^{\prime}(\omega_c)\bigg\} \vspace{1ex}\\
&=&sign\bigg\{2\omega_c\bigg(-C^{\prime}(\omega_c)
-\delta\bigg(\displaystyle\frac{S(\omega_c)}{\omega_c}\bigg)^{\prime}\bigg)\bigg\} \vspace{1ex}\\
&=&sign\bigg\{-C^{\prime}(\omega_c)-\delta\bigg(\displaystyle\frac{S(\omega_c)}
{\omega_c}\bigg)^{\prime}\bigg\}.
\end{array}
\end{displaymath}
This concludes the proof. \qquad
\end{proof}

\emph{Remark 2.} Consider the function $g$ defined by
(\ref{equationg}) and denote by $\alpha$ the quantity
\begin{displaymath}
\alpha:=\bigg(2\displaystyle\int_0^{\overline{\tau}}
e^{-\gamma\tau}f(\tau)d\tau-1 \bigg)\beta(0).
\end{displaymath}
Define the sets
\begin{displaymath}
\Omega:=\{ \omega>0 ; \ 0<g(\omega)<\alpha \ \textrm{ and } \
g^{\prime}(\omega)=0 \} \qquad \textrm{ and } \qquad
\Lambda:=g(\Omega).
\end{displaymath}
One can notice that $\Lambda$ is finite (or empty). If
$\delta\in(0,\alpha)\setminus\Lambda$, then
\begin{displaymath}
\frac{d Re(\lambda)}{d(-\beta^*)}\bigg|_{\beta^*=\beta^*_c}\neq0.
\end{displaymath}
Indeed, we have
\begin{displaymath}
g^{\prime}(\omega)=
-\frac{\omega}{S(\omega)}\bigg(g(\omega)\bigg(\frac{S(\omega)}
{\omega}\bigg)^{\prime}+C^{\prime}(\omega)\bigg), \qquad \omega>0.
\end{displaymath}
Since $\delta\notin\Lambda$, we have $g^{\prime}(\omega_c)\neq0$.
Moreover, $g(\omega_c)=\delta$. Thus
\begin{displaymath}
C^{\prime}(\omega_c)\neq-\delta\bigg(\displaystyle\frac{S(\omega_c)}{\omega_c}\bigg)^{\prime}.
\end{displaymath}
We conclude by using (\ref{conditionhopfbifurcation}).

Lemma \ref{lemmaexistenceproperties}, together with Remark 2,
allows us to state and prove the following theorem.

\begin{theorem} \label{theoremhopfbifurcation}
Assume that the function $\tau\mapsto e^{-\gamma\tau}f(\tau)$ is
decreasing. Then, for each $\delta\notin\Lambda$ satisfying
(\ref{condexistenceequilibrium}), there exists
$\beta^*_c<\widetilde{\delta}$ such that the equilibrium $x\equiv
x^*$ is locally asymptotically stable when
$\beta^*_c<\beta^*\leq\widetilde{\delta}$ and a Hopf bifurcation
occurs at $x\equiv x^*$ when $\beta^*=\beta^*_c$.
\end{theorem}

\begin{proof}
First, recall that $x\equiv x^*$ is locally asymptotically stable
when $\beta^*=\widetilde{\delta}$ (see Theorem
\ref{theoremstability}). We recall that, from the properties of
the function $g$, equation (\ref{equationg}) has a finite number
of solutions (see Lemma \ref{lemmaexistenceproperties}). We set
\begin{displaymath}
\beta^*_c = -\frac{\delta}{1-2C(\omega^*_c)},
\end{displaymath}
where $\omega^*_c$ is the smaller positive real such that
\begin{displaymath}
C(\omega^*_c)=\min\{C(\omega); \ \omega \textrm{ is a solution of
} (\ref{equationg})\}.
\end{displaymath}
Then, $\beta^*_c$ is the maximum value of $\beta^*$ (as defined in
Lemma \ref{lemmaexistenceproperties}) which gives a solution of
(\ref{systempureimaginaryroot}). From Lemma
\ref{lemmaexistenceproperties}, (\ref{characteristicequation}) has
no purely imaginary roots while
$\beta^*_c<\beta^*\leq\widetilde{\delta}$. Consequently,
Rouch\'e's Theorem [\ref{rouche}, p.248] leads to the local
asymptotic stability of $x\equiv x^*$.

When $\beta^*=\beta^*_c$, (\ref{characteristicequation}) has a
pair of purely imaginary roots $\pm i\omega_c$, $\omega_c>0$ (see
Lemma \ref{lemmaexistenceproperties}). Moreover, since
$\delta\notin\Lambda$, Remark 2 implies that
\begin{displaymath}
\frac{d Re(\lambda)}{d(-\beta^*)}\bigg|_{\beta^*=\beta^*_c}\neq0.
\end{displaymath}
Assume, by contradiction, that
\begin{displaymath}
\frac{d Re(\lambda)}{d(-\beta^*)}<0
\end{displaymath}
for $\beta^*>\beta^*_c$, $\beta^*$ close to $\beta^*_c$. Then
there exists a characteristic root $\lambda(-\beta^*)$ such that
$Re \lambda(-\beta^*)>0$. This contradicts the fact that $x\equiv
x^*$ is locally asymptotically stable when $\beta^*>\beta^*_c$.
Thus, we obtain
\begin{displaymath}
\frac{d Re(\lambda)}{d(-\beta^*)}\bigg|_{\beta^*=\beta^*_c}>0.
\end{displaymath}
This implies the existence of a Hopf bifurcation  at $x\equiv x^*$
for $\beta^*=\beta^*_c$. \qquad
\end{proof}

\begin{figure}[!hpt]
\begin{center}
\includegraphics[width=10cm, height=8cm]{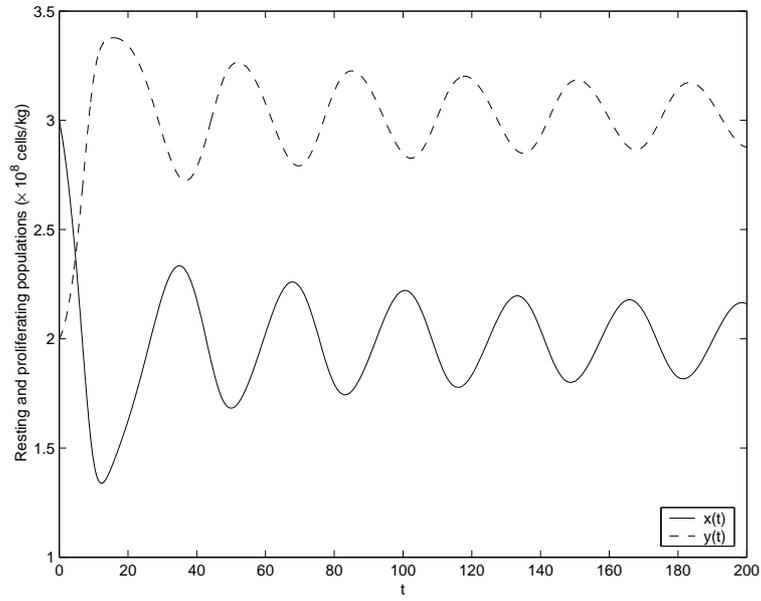}
\end{center}
\caption{The solutions of system (\ref{problem}), $x(t)$ (solid
curve) and $y(t)$ (dashed curve), are drawn when the Hopf
bifurcation occurs. This corresponds to $n=2.53$ with the other
parameters given by (\ref{parameters}) and $\overline{\tau}=7$
days. Periodic solutions appear with period of the oscillations
about $33$ days.}\label{hopfbiffig}
\end{figure}

\begin{figure}[!hpt]
\begin{center}
\includegraphics[width=10cm,height=8cm]{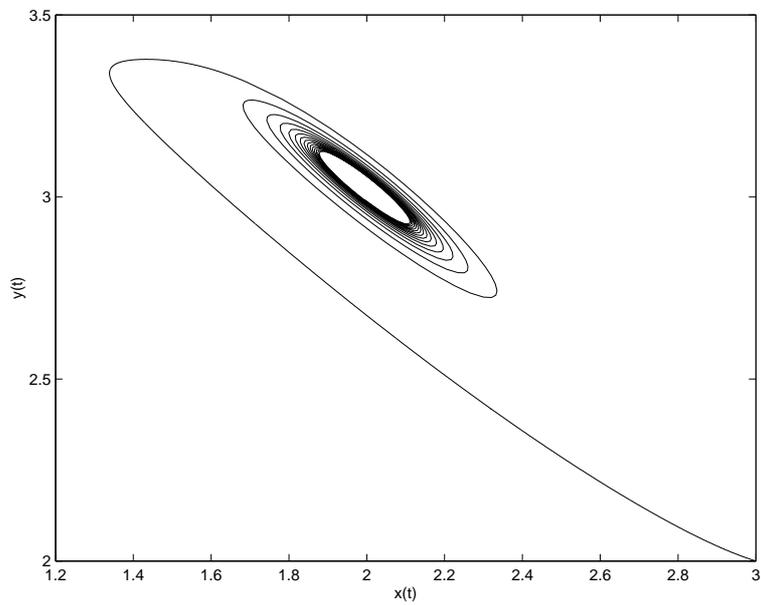}
\end{center}
\caption{For the values used in Fig.\ref{hopfbiffig}, the
solutions are shown in the $(x,y)$-plane: the trajectories reach a
limit cycle, surrounding the equilibrium.}\label{hopfbiffig2}
\end{figure}

With the values of $\delta$, $\gamma$ and $\beta_0$ given by
(\ref{parameters}), and $\overline{\tau}=7$ days, (\ref{problem})
has periodic solutions for $\beta^*_c=-0.3881$ with a period about
$33$ days. This value of $\beta^*_c$ corresponds to $n=2.53$ (see
Fig. \ref{hopfbiffig} and \ref{hopfbiffig2}). The function $f$ is
given by (\ref{uniformlaw}).


The bifurcation parameter was chosen to be $\beta^*$ in this
study, and the values of $\beta^*$ depend strongly on the
sensitivity $n$ of the function $\beta(x)$, since all other
parameters are fixed by (\ref{condexistenceequilibrium}). In this
model, the sensitivity $n$ plays a crucial role in the appearance
of periodic solutions. Pujo-Menjouet and Mackey \cite{pujomackey}
already noticed the influence of this parameter on system
(\ref{problem}) when the delay is constant (or equivalently, when
$f$ is a Dirac measure). The sensitivity $n$ describes the way the
rate of introduction in the proliferating phase reacts to changes
in the resting phase population produced by external stimuli: a
release of erythropoietin, for example, or the action of some
growth factors.

Of course, the influence of other parameters (like mortality rates
$\delta$ and $\gamma$, or the minimum and maximum delays
$\underline{\tau}$ and $\overline{\tau}$) on the appearance of
periodic solutions could be studied. However, since periodic
hematological diseases --- defined and described in Section
\ref{sectiondiscussion} --- are supposed to be due to hormonal
control destabilization (see \cite{fm99}), then the parameter $n$,
among other parameters, seems to be appropriate to identify causes
leading to periodic solutions in (\ref{problem}).

\section{Discussion}\label{sectiondiscussion}

Among the wide range of diseases affecting blood cells, periodic
hematological diseases (Haurie {\it et al}. \cite{hdm98}) are of
main importance because of their intrinsic nature. These diseases
are characterized by significant oscillations in the number of
circulating cells, with periods ranging from weeks (19--21 days
for cyclical neutropenia \cite{hdm98}) to months (30 to 100 days
for chronic myelogenous leukemia \cite{hdm98}) and amplitudes
varying from normal to low levels or normal to high levels,
depending on the cells types \cite{hdm98}. Because of their
dynamic character, periodic hematological diseases offer an
opportunity to understand some of the regulating processes
involved in the production of hematopoietic cells, which are still
not well-understood.

Some periodic hematological diseases involve only one type of
blood cells, for example, red blood cells in periodic auto-immune
hemolytic anemia (B\'elair {\it et al.} \cite{bmm95}) or platelets
in cyclical thrombocytopenia (Santillan {\it et al.}
\cite{sbmm00}). In these cases, periods of the oscillations are
usually between two and four times the bone marrow production
delay. However, other periodic hematological diseases, such as
cyclical neutropenia (Haurie {\it et al}. \cite{hdm98}) or chronic
myelogenous leukemia (Fortin and Mackey \cite{fm99}), show
oscillations in all of the circulating blood cells, i.e., white
cells, red blood cells and platelets. These diseases involve
oscillations with quite long periods (on the order of weeks to
months). A destabilization of the pluripotential stem cell
population (from which all of the mature blood cells types are
derived) seems to be at the origin of these diseases.

We focus, in particularly, on chronic myelogenous leukemia (CML),
a cancer of the white cells, resulting from the malignant
transformation of a single pluripotential stem cell in the bone
marrow (Pujo-Menjouet {\it et al.} \cite{mackeypujo2}). As
described in Morley {\it et al.} \cite{morley}, oscillations can
be observed in patients with CML, with the same period for white
cells, red blood cells and platelets. This is called periodic
chronic myelogenous leukemia (PCML). The period of the
oscillations in PCML ranges from 30 to 100 days (Haurie {\it et
al}. \cite{hdm98}, Fortin and Mackey \cite{fm99}), depending on
patients. The difference between these periods and the average
pluripotential cell cycle duration (between 1 to 4 days, as
observed in mouses, see Mackey \cite{mackey2001}) is still not
well-understood.

Recently, to understand the dynamics of periodic chronic
myelogenous leukemia, Pujo-Menjouet {\it et al.}
\cite{mackeypujo2} considered a model for the regulation of stem
cell dynamics and investigated the influence of parameters in this
stem cell model on the oscillations period when the model becomes
unstable and starts to oscillate. In this paper, taking into
account the fact that a cell cycle has two phases, that is, stem
cells in process are either in a resting phase or actively
proliferating, and assuming that cells divide at different ages,
we proposed a system of differential equations with distributed
delay to model the dynamics of hematopoietic stem cells. By
constructing a Lyapunov functional, we gave conditions for the
trivial equilibrium to be globally asymptotically stable. Local
stability and Hopf bifurcation of the nontrivial equilibrium were
studied, the existence of a Hopf bifurcation leading to the
appearance of periodic solutions in this model, with a period
around 30 days at the bifurcation.

\begin{figure}[!hpt]
\begin{center}
\includegraphics[width=10cm, height=8cm]{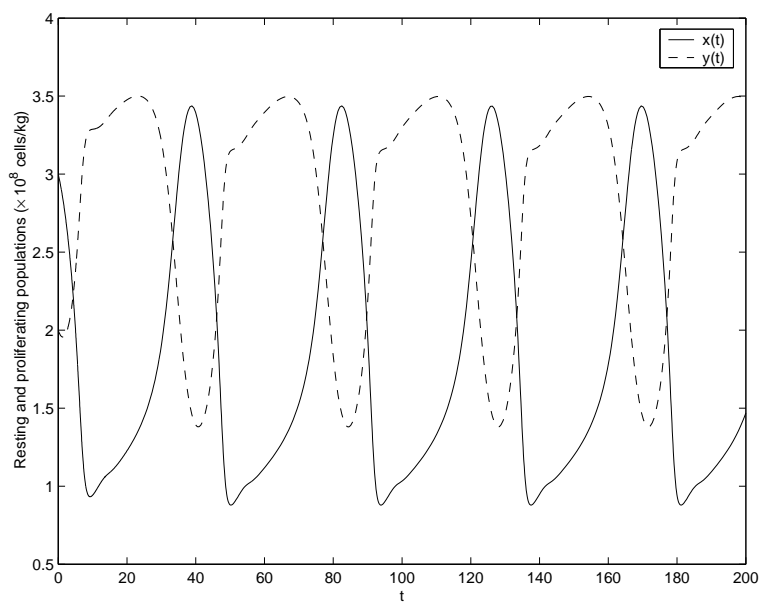}
\end{center}
\caption{Solutions $x(t)$ (solid curve) and $y(t)$ (dashed curve)
of system (\ref{problem}) oscillate with periods close to 45 days;
the parameters are the same as in Fig \ref{hopfbiffig}, with
$n=3$. The amplitudes of the oscillations range from low values to
normal values.}\label{figoscillations1}
\end{figure}

Numerical simulations show that periodic solutions occur after the
bifurcation, with periods increasing as the bifurcation parameter
(the sensitivity $n$) increases. In Fig. \ref{figoscillations1},
solutions oscillate around the equilibrium values with periods
around 45 days. Moreover, amplitudes of the oscillations range
from low values to normal values. The sensitivity is equal to
$n=3$, that is the parameters are given by (\ref{parameters}).
This corresponds to values given by Mackey \cite{mackey1978},
values for which abnormal behavior (periodic) is usually observed
in all circulating blood cells types.

\begin{figure}[!hpt]
\begin{center}
\includegraphics[width=10cm, height=8cm]{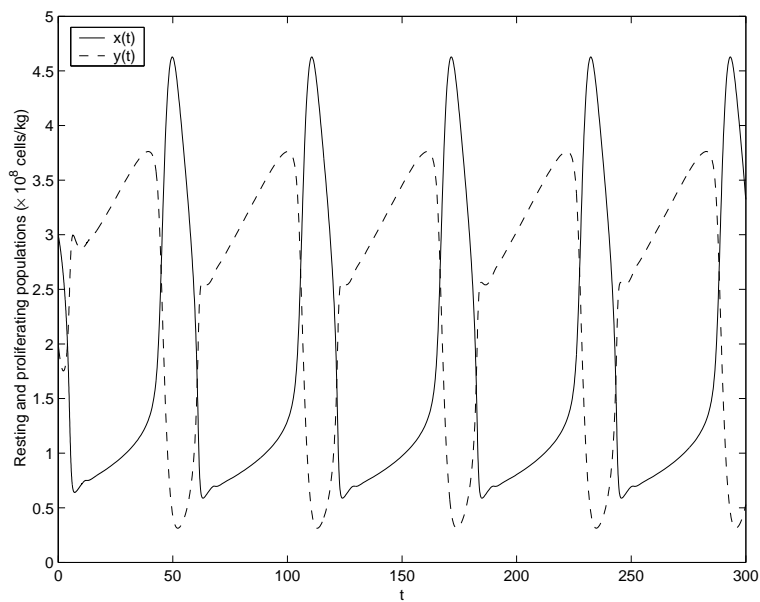}
\end{center}
\caption{Solutions $x(t)$ (solid curve) and $y(t)$ (dashed curve)
of system (\ref{problem}) oscillate with periods close to 70 days;
the parameters are the same as in Fig \ref{hopfbiffig}, with
$n=4$. The amplitudes of the oscillations range from low values to
high values.}\label{figoscillations2}
\end{figure}

When $n$ continues to increase, longer oscillations periods are
observed with amplitudes varying from low values to high values
(see Fig. \ref{figoscillations2}). This situation characterizes
periodic chronic myelogenous leukemia, with periods in the order
of two months (70 days).

Moreover, the oscillations observed in Fig. \ref{figoscillations1}
and \ref{figoscillations2} look very much like relaxation
oscillations. Experimental data from patients with PCML suggest
that the shape of oscillations is of a relaxation oscillator type
\cite{fm99, hdm98}. Furthermore, Fowler and Mackey \cite{fm02}
showed that a model for hematopoiesis with a discrete delay may
also exhibit relaxation oscillations. Therefore, it seems that not
only periods and amplitudes of the oscillations correspond to the
ones observed in PCML, but also the shape of the oscillations.

Numerical simulations demonstrated that long period oscillations
in the circulating cells are possible in our model even with short
duration cell cycles. Thus, we are able to characterize some
hematological diseases, especially those which exhibit a periodic
behavior of all the circulating blood cells.

\end{document}